\documentclass[preprint,12pt]{elsarticle}

\usepackage{graphicx}
\usepackage{caption}
\usepackage{subcaption}
\usepackage{comment}
\usepackage{algorithm}
\usepackage{algpseudocode}

\usepackage{amsmath} % assumes amsmath package installed
\usepackage{amssymb}  % assumes amsmath package installed
\usepackage{amsthm}  % assumes amsmath package installed
\usepackage{bm}
\usepackage{enumitem}

\usepackage{lipsum}
\usepackage{color}

\newtheorem{prop}{Proposition}
\newtheorem{defn}{Definition}
\newtheorem{cor}{Corollary}

\newtheorem{thm}{Theorem}

\numberwithin{equation}{section}

\newcommand{\R}{\mathbb{R}}

\newcommand{\de}{\mathrm{d}}
\newcommand{\x}{x}

\DeclareMathOperator*{\argmin}{arg\,min}

\usepackage{amssymb}
\journal{}

\begin{document}

\begin{frontmatter}



\title{Hamilton--Jacobi Based Policy-Iteration via Deep Operator Learning}


\author[inst1]{Jae Yong Lee}

\affiliation[inst1]{organization={Department of Artificial Intelligence, Chung-Ang University},%Department and Organization
            addressline={84 Heukseok-ro}, 
            city={Dongjak-gu},
            postcode={06974}, 
            state={Seoul},
            country={Republic of Korea}}

\author[inst2]{Yeoneung Kim}

\affiliation[inst2]{organization={Department of Applied Artificial Intelligence, SeoulTech},%Department and Organization
            addressline={232 Gongneung-ro}, 
            city={Nowon-gu},
            postcode={01811}, 
            state={Seoul},
            country={Republic of Korea}}

\begin{abstract}
The framework of deep operator network (DeepONet) has been widely exploited thanks to its capability of solving high dimensional partial differential equations. In this paper, we incorporate DeepONet with a recently developed policy iteration scheme to numerically solve optimal control problems and the corresponding Hamilton--Jacobi--Bellman (HJB) equations. A notable feature of our approach is that once the neural network is trained, the solution to the optimal control problem and HJB equations with different terminal functions can be inferred quickly thanks to the unique feature of operator learning. Furthermore, a quantitative analysis of the accuracy of the algorithm is carried out via comparison principles of viscosity solutions. The effectiveness of the method is verified with various examples, including 10-dimensional linear quadratic regulator problems (LQRs).
\end{abstract}


% \begin{highlights}
% \item We incorporate the concept of a physics-informed deep operator network (PI-DeepONet) with a novel policy iteration scheme to effectively solve optimal control problems with various terminal functions. 
% \item The error analysis is provided through the comparison principle for viscosity solutions.
% \item The effectiveness of our new algorithm is verified through various numerical experiments, including a vehicle problem and high-dimensional linear quadratic regulator (LQR) problems with a compact control set.
% \end{highlights}


\end{frontmatter}



\section{Introduction}
Optimal control has been an important area of study in many fields due to its interesting mathematical properties and applicability for solving real-world problems, including portfolio optimization, economics, robot control, operational research, and the design of optimal vaccination strategies~\cite{lenhart2007optimal,lin2007adaptive,sandholm2021hamilton,mitchell2011scalable,barles1997convergence,forsyth2007numerical}. Over the years, there have been numerous attempts to understand optimal control problems from various perspectives, among which the Pontryagin Maximum Principle and Bellman's dynamic programming principle (DPP)~\cite{bellman1966dynamic} are fundamental. In particular, in the regime of continuous time and space, the value function of optimal control problems leads to the Hamilton--Jacobi--Bellman (HJB) partial differential equation that is known to have the unique viscosity solution~\cite{crandall1992user}. Once the solution to the HJB is known, it is leveraged to synthesize an optimal control if the value function is smooth~\cite{bardi1997optimal,evans2022partial}. However, the explicit representation of solutions to HJB may not exist, so it is essential to consider numerical approximation.

To understand the viscosity solution to the HJB equations as well as corresponding optimal controls numerically, various methods have been proposed via discretization of the state space~(e.g., \cite{osher1991high, jiang2000weighted, barles2002convergence,forsyth2007numerical}). Unfortunately, the convergence of numerical methods often suffers from scalability problems as the computational complexity increases exponentially with the dimension of the state space, which is referred to as \textit{the curse of dimensionality}. Many approaches for approximating the value function and optimal control have been proposed to avoid such a phenomenon. In particular, when the Hamiltonian depends only on the costate, we have a simple representation of the solution, called the Hopf--Lax formula~\cite{hopf1965generalized,bardi1984hopf}. Various numerical algorithms based on the representation formulas are shown to be effective for some classes of high-dimensional HJB equations~\cite{darbon2016algorithms, aliyu2016modified,chow2019algorithm, darbon2020decomposition, lee2020hopf, yegorov2021perspectives,kim2022representation}. 
Additionally, \cite{chen2021lax} tries to identify explicit solutions to a certain class of optimal control problems using the Hopf--Lax formula. Another notable work is \cite{lee2021computationally}, where a generalized Lax formula is obtained via DPP to handle nonlinear problems. Unfortunately, control trajectories obtained by this method often present rapid oscillation. As an alternative, neural network-based approaches for solving high dimensional HJB equations and synthesizing optimal control problems have gained lots of attention. Pioneered by~\cite{han2018solving} where the authors implement deep neural networks to solve high-dimensional partial differential equations including Hamilton--Jacobi (HJ) equations, some neural network architectures are proposed to solve some state-independent high-dimensional HJ equations numerically without discretization~\cite{MR3874585,darbon2021some,darbon2020overcoming}. 

Recently, a new research area called operator learning has emerged, focusing on learning mappings from functions to functions to solve partial differential equations~\cite{lu2019deeponet,lu2021learning,li2020fourier}. Operator learning aims to acquire the mathematical operators that govern the behavior of a physical system from data. These methods rely on a pre-generated dataset of input-output pairs, obtained by solving the PDEs numerically. Once the model is trained using this dataset, it enables real-time predictions for diverse input data.

However, generating a reliable training dataset for these methods requires a substantial number of numerical solutions obtained through extensive computations. To address this challenge, the framework of physics--informed machine learning methods has emerged, enabling the learning of operators without relying on pre-generated data. Notably, \cite{PIDON} introduced an extended version of the deep operator network (DeepONet) \cite{lu2021learning}, called the physics-informed DeepONet (PI-DeepONet). This model is enhanced by incorporating the residual of the PDE as an additional loss function.

In this paper, a novel way of solving deterministic optimal control problems and HJB equations is proposed by combining the framework of PI-DeepONet with the policy iteration scheme developed by~\cite{tang2023policy}. In~\cite{tang2023policy}, the authors propose an idea of iterating the unique viscosity solutions to linearized HJB equations and show that the scheme guarantees the convergence to the unique solution of the original HJB equation. Motivated by this, we introduce the DeepONet structure that recursively solves the linearized HJB equations without discretizing the state space. Another great advantage of our method is that the DeepONet architecture can be trained to infer the solution to HJB with different terminal functions. Some quantitative results on the accuracy of approximate solutions are provided, and 
we verify the effectiveness of our algorithm in solving HJB equations and synthesizing optimal control in various dimensions. 

%The main idea is to incorporate DeepONet with recursively defined linear partial differential equations to efficiently solve Hamilton--Jacobi equations that can infer not only optimal controls but also viscosity solutions to HJB for multiple terminal value functions at once. 

\subsection{Organization of paper}
The paper is organized as follows. In Section~\ref{sec:pre}, we provide a gentle introduction to the policy iteration scheme and the framework of PI-DeepONet. In the following section, we introduce our new algorithm and present the error analysis. In Section~\ref{sec:exp}, we present numerical results to support the algorithm. Finally, we conclude in Section~\ref{sec:conc}.

\subsection{Notations}
Let us introduce notations used throughout the paper
\begin{itemize}
\item For $d \in\mathbb {N}$, $x \in \mathbb{R}^d$ and $p \geq 1$, we define 
\[
\|x\|_p :=  \sum_{i=1}^n |x[i]|^p)^{1/p}, 
\]
where $i$th component of $x$ is denoted by $x[i]$ ($i=1.,,,.d$).
\item Let $C^1(X)$ denote continuously differentiable functions defined in $X$.

\item For $f\in C^1(\mathbb{R}^d)$, 
\[
D_x f(x):= (\frac{\partial f}{\partial x[1]},...,\frac{\partial f}{\partial x[d]})
\quad
\text{and}\quad
\Delta f := \sum_{i=1}^d \frac{\partial^2 f}{\partial x[i]^2}.
\]
Given $h>0$, we define
%\[
%\nabla^h f(x):=\bigg(\frac{f(x+he_1)-f(x-he_1)}{2h},...,\frac{f(x+h e_d)-f(x-h e_d)}{2h}\bigg)
%\]
\[
\nabla^h f(x):=\bigg(\nabla^h_1 f (x),...,\nabla^h_d f (x)\bigg)
\]
for
\[
\nabla^h_i f (x) := \frac{f(x+h e_i)-f(x-h e_i)}{2h}
\]
and
\[
\Delta^h f(x):= \sum_{i=1}^d \frac{f(x+he_i) -2f(x)+f(x-he_i)}{h^2}.
\]
\end{itemize}

\section{Preliminary}\label{sec:pre}
\subsection{Optimal control and Hamilton--Jacobi--Bellman equation}
Consider a continuous-time deterministic dynamical system of the form
\begin{equation*}
\begin{split}
\frac{dx_t}{dt} &= f(t, x_t, u_t), \quad   t \in (0,T),
\end{split}
\end{equation*}
where $x_t \in \R^d$ and $u_t \in U \subset \R^m$ are the system state and control at time $t$ respectively. The set of admissible controls is defined as
\[
\mathcal{U}_t:=\{u:[t,T] \rightarrow U  :  \text{$u$ is measurable}\},
\]
where the control set $U$ is assumed to be compact.

%A finite-horizon optimal control problem can then be formulated as 
%\[
%\min_{u \in \mathcal{U}} \left \{ \int_0^T L(t,x(t),u(t)) \de t+g(x(T)) : x(0) = x \right \},
%\]
%where $L$ is the stage-wise cost function, and $g$ is the terminal cost function. 
%Throughout the paper, we always assume that

Let $V(t,x)$ denote the optimal value function defined as
\begin{equation*}
V(t,x):= \inf_{u \in \mathcal{U}_t} J(t,x,u),
\end{equation*}
where
\[
J(t,x,u):=   \int_t^T L(s,x_s,u_s) \de s+g(x_T) \quad \text{given}\quad x_t = x.
\]
The value function represents the optimal cost-to-go starting from $x_t = x$ at time $t$. 
Under the assumption that $f$, $L$, and $g$ are uniformly Lipscthiz continuous, it is well known that $V$ is the unique viscosity solution of the following HJ equation~\cite{tran2021hamilton}:
\begin{equation}\label{eq:hj}
\begin{cases}
   \partial_t V + H(t, x,D_x V) = 0 &\quad\text{in}\quad(0,T)\times\R^d,\\
    V(T,x)=g(x)&\quad\text{on}\quad\R^d,
\end{cases}
\end{equation}
where the Hamiltonian is given by
\[
H(t,x,p):=\inf_{{\boldsymbol{u}} \in U} \left\{ p\cdot f(t,x, {\boldsymbol{u}})+L(t, x, {\boldsymbol{u}}) \right\}.
\]
If a solution $V$ to~\eqref{eq:hj} is continuously differentiable, the optimal policy is given by
\[
u_*(t,x)=u(t,x,D_x V(t,x)),
\]
where
\begin{equation}\label{eq:opt}
u(t,x,p):=\argmin_{\boldsymbol{u}\in U} \{p \cdot f(t,x,{\boldsymbol{u}})+L(t,x,{\boldsymbol{u}})\}.
\end{equation}
For verification of this property, we refer~\cite{evans2022partial, bardi1997optimal,tran2021hamilton} and references therein. 

Throughout the paper, we assume that $p\cdot f(t,x,\boldsymbol{u})+L(t,x,\boldsymbol{u})$ is convex in $\boldsymbol{u}$ for each $(t,x,p)$ and the optimization problem~\eqref{eq:opt} can be solved using standard convex optimization tools. The assumption seems restrictive but it still covers various classes of optimal control problems such as control affine systems where the dynamics and Lagrangian are given by 
\[
f(t,x,\boldsymbol{u})=g_1(t,x)+g_2(t,x)\cdot\boldsymbol{u} \quad\text{and}\quad L(x,\boldsymbol{u})=\|x\|_2^2+\|\boldsymbol{u}\|_2^2
\]
for some $g_1,g_2 \in C([0,T]\times \mathbb{R}^d)$.

\subsection{Semi-discrete iterative scheme}
We revisit the semi-discrete scheme for the policy iteration proposed by~\cite{tang2023policy} where the following assumptions are necessary for the well-posedness of linearized equations.
\begin{enumerate}[label=(A\arabic*)]
\item $L(\cdot,\cdot,\cdot)$, $f(\cdot,\cdot,\cdot)$ and $g(\cdot)$ are uniformly bounded and Lipschitz continuous;
\item $u_0:\R \times \R^d \rightarrow U$ and $u$ defined in~\eqref{eq:opt} are Lipschitz continuous in all arguments. Furthermore,~\eqref{eq:opt} achieves a unique minimum for each $t$ and $x$.
\end{enumerate}

Let $u_0$ be a continuous function satisfying (A2). For $h\in (0,1)$, $N>0$ and $n \in \mathbb{N}$, we consider
\begin{equation}\label{eq:semidiscret}
    \begin{cases}
    \partial_t V_n^h + L(t,x,u_n)+\nabla^h V_n^h \cdot f(t,x,u_n)=-Nh \Delta^h V_n^h \quad&\text{in}\quad (0,T)\times \mathbb{R}^d,\\
    V_n^h(T,x)=g(x) \quad&\text{on}\quad \mathbb{R}^d.
\end{cases}
\end{equation}
where
\begin{equation}\label{eq:update}
u_{n+1}(t,x)=u(t,x,\nabla^h V_n^h)=\argmin_{{\boldsymbol{u}}\in U} \{\nabla^h V_n^h \cdot f(t,x,{\boldsymbol{u}})+L(t,x,{\boldsymbol{u}})\}.
\end{equation}
It is known that the sequence of solutions is monotone, that is, $V_n^{h} \geq V_{n+1}^h$ for all $n$. Furthermore, starting from arbitrary $u_0$, the following convergence property of $V_n^h$ is well-established in~\cite{tang2023policy}.
\begin{thm}
Let $N \geq \max\{1,\|f\|_\infty/2\}$ and $h \in (0,1)$ be given and assume that (A1) and (A2) are enforced. For $V_n^h$ solving~\eqref{eq:semidiscret}, we have that
\[
V_n^h(t,x) \rightarrow V^h(t,x) \quad\text{locally uniformly}\quad\text{as}\quad n\rightarrow \infty,
\]
where $V^h$ is a viscosity solution to
\begin{equation}\label{eq:app_hj}
    \begin{cases}
    \partial_t V^h +H(t,x,\nabla^h V^h)=-N h \Delta^h V^h\quad&\text{in}\quad (0,T)\times \mathbb{R}^d,\\
    V^h(T,x)=g(x)\quad&\text{on}\quad \mathbb{R}^d.
    \end{cases}
\end{equation}
Here, $H(t, x, p):=\inf_{{\boldsymbol{u}}\in U} \left\{ p\cdot f(t,x, {\boldsymbol{u}})+L(t, x, {\boldsymbol{u}})\right\}$. 
 
Furthermore, for the unique viscosity solution $V$ of~\eqref{eq:hj}, we deduce that
\[
|V-V^h| \leq C \sqrt{h}\quad \text{in} \quad [0,T]\times \mathbb{R}^d
\]
for a global constant $C>0$.
\end{thm}
For details, we refer readers to~Section 3 of~\cite{tang2023policy}. A salient feature of the policy iteration scheme is that one can simply alternate between linear semi-discrete partial differential equation~\eqref{eq:semidiscret} and policy update~\eqref{eq:update} to approximate the solution to the original HJB equation~\eqref{eq:app_hj}. We now introduce the framework of PI-DeepONet in the following subsection.
%Motivated by PI, we incorporate a physics-informed operator learning framework with PI to develop a method for solving HJ equation for any terminal value function $g(x)$ given. 

\subsection{Physics-Informed DeepONet (PI-DeepONet)}
\begin{figure}[h]
    \centering
    \includegraphics[width=\textwidth]{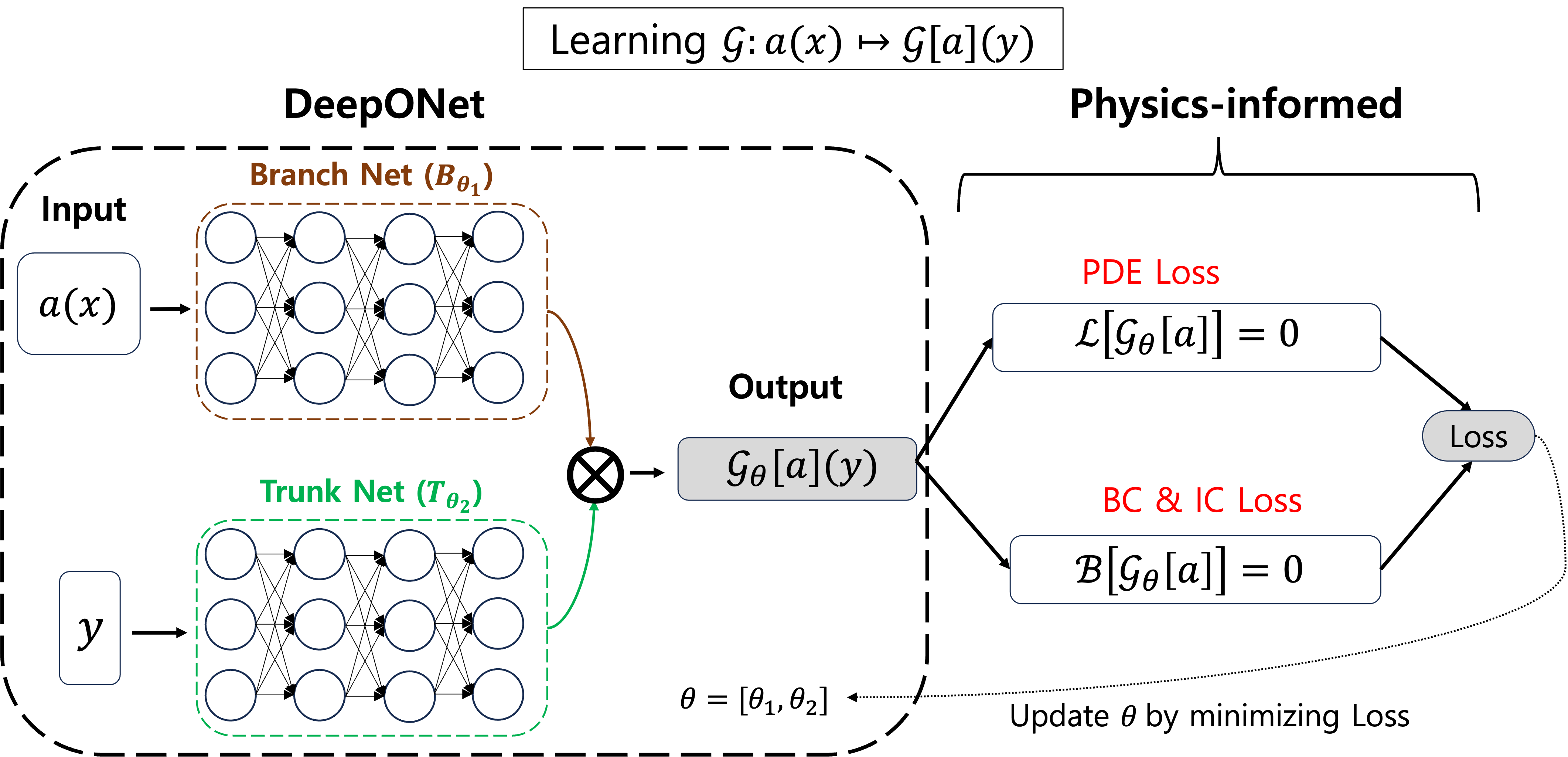}
    \caption{PI-DeepONet framework.}
    \label{fig:pideeponet}
\end{figure}
Lu et al. \cite{lu2019deeponet} proposed a DeepONet architecture based on the universal approximation theorem for operators~\cite{chen1995universal}. 
The left square in Figure \ref{fig:pideeponet} describes the DeepONet structure to learn the operator $\mathcal{G}:\mathcal{A}\to \mathcal{F}$ such that $\mathcal{G}:a(x)\mapsto \mathcal{G}[a](y)$ where $\mathcal{A}$ and $\mathcal{F}$ are Banach spaces.
%Altogether,  Here, $\mathcal{G}[a](y)\in\mathcal{F}$ is the solution of PDE $\mathcal{L}[\mathcal{G}[a]]=0$ with the boundary condition $\mathcal{B}[\mathcal{G}[a]]=0$. 
The DeepONet $\mathcal{G}_\theta$, which approximates the true operator $G$, consists of two neural networks: a branch net and a trunk net parametrized by $\theta_1$ and $\theta_2$ respectively. The branch net takes an input $\vec{a}=[a(x_1),a(x_2),...,a(x_k)]$ at fixed sensor points $x_1,x_2,...,x_k$ and generates $p$-values $B_{\theta_1}(\vec{a})=[B_{\theta_1}(\vec{a})_1,...,B_{\theta_1}(\vec{a})_p]$ as an output. The trunk net takes a variable $y$ for the target function space and generates $p$-values $T_{\theta_2}(y)=[T_{\theta_2}(y)_1,...,T_{\theta_2}(y)_p]$ as an output. The DeepONet approximates the target function $\mathcal{G}[a](y)$ by taking the inner product of the outputs of these two networks, that is,
\begin{equation}
    \mathcal{G}[a](y)\approx \sum_{i=1}^p B_{\theta_1}(\vec{a})_i T_{\theta_2}(y)_i.
\end{equation}

To train DeepONet, it is essential to have input-output data pairs. Learning DeepONet presents significant limitations when data is unavailable or obtaining data incurs substantial computational costs. To overcome this limitation, PI-DeepONet \cite{PIDON} has emerged, facilitating learning without the need for data. The right side of Figure \ref{fig:pideeponet} illustrates the physics-informed loss used in the PI-DeepONet structure to train DeepONet without input-output data pairs. Specifically, we can identify solving boundary value PDE as finding operator $\mathcal{G}_\theta$ such that 
\begin{equation}
\begin{cases}
\mathcal{L}[\mathcal{G}_{\theta}[a]]=0, \quad \quad &\text{(PDE operator)},\\
\mathcal{B}[\mathcal{G}_{\theta}[a]]=a \in\mathcal{A},\quad \quad&\text{(Boundary operator).}
\end{cases}
\end{equation}
The learning of $\mathcal{G}_\theta$ proceeds by calculating the PDE loss which is the residual of the equation arising from the approximate solution $\mathcal{G}_\theta[a](y)$, and the boundary condition loss via automatic differentiation tool \cite{MR3800512,paszke2017automatic}. 

% To make use of the monotonicity property of the semi-discrete scheme,
% \[
% v_{h+1}(t,x) \leq v_h (t,x)\quad\text{in}\quad(0,T)\times\mathbb{R}^d,
% \]
% we add 
\section{Main Results and Methodology}
Let $h \in (0,1)$ in the rest of the paper. We also assume that $N \geq \max\{1,\|f\|_\infty/2$\}, $L$ and $f$ satisfying Assumption (A1) are given. Inspired by the iterative scheme~\eqref{eq:semidiscret}, we consider an operator $\mathcal{H}$ that maps $g$ to $V^h$ solving 
\begin{equation*}
    \begin{cases}
    \partial_t V^h + H(t,x,\nabla^h V^h)=-Nh \Delta^h V^h \quad\text{in}\quad (0,T)\times \mathbb{R}^d,\\
    V^h(T,x)=g(x)\quad\text{on}\quad\mathbb{R}^d
\end{cases}
\end{equation*}
for any $g$. The condition $N\geq \|f\|_\infty/2$ is necessary to guarantee the monotonicity of the numerical Hamiltonian as found in~\cite{bardi1997optimal,achdou2013hamilton,crandaltwoapp}. 

From this perspective, we aim at finding an operator $\mathcal{H}_\theta$ with the learnable parameter $\theta$ to approximate $\mathcal{H}$ for any input function $g$, that is,
\[
V^h=\mathcal{H}[g] \approx \mathcal{H}_\theta[g].
\]

% %that approximates the operator $\mathcal{G}[g]=V^h$.
As a result, we deduce that $\mathcal{H}_\theta[g]$ approximates $V$ solving~\eqref{eq:hj} since
\[
|V - \mathcal{H}_\theta[g] | \leq |V - V^h| +  |V^h - \mathcal{H}_\theta[g]|.
\]

In the following subsection, we describe the learning procedure by iteratively applying an operator. It is also verified that small training errors would yield an accurate solution, which comes from simple error analysis based on the comparison principle established in the theory of viscosity solutions.

\subsection{Solving optimal control by using PI-DeepONet iteratively}

\begin{figure}[h]
    \centering
    \includegraphics[width=\textwidth]{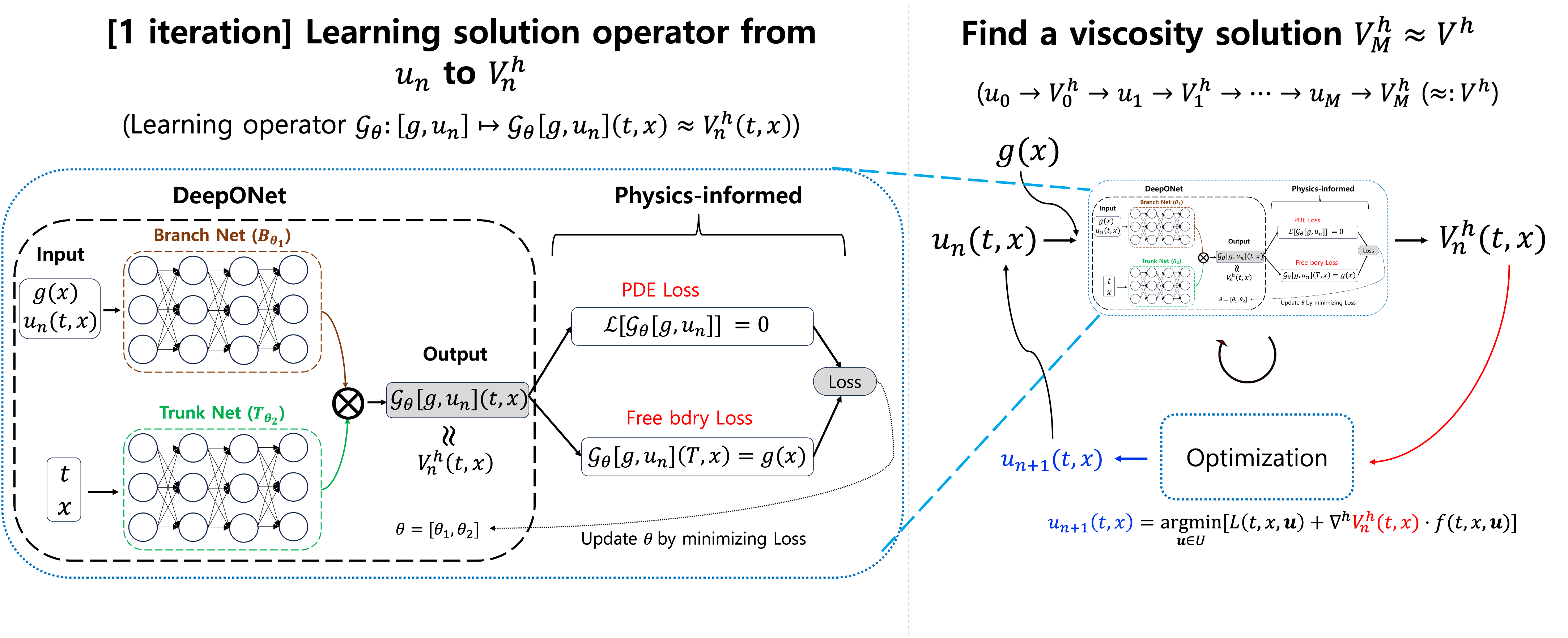}
    \caption{Our framework: Hamilton-Jacobi based policy-iteration via PI-DeepONet}
    \label{fig:framework}
\end{figure}

\begin{algorithm}[h]
\caption{Iterative PI-DeepONet training procedure}
\label{alg:main}
\begin{algorithmic}[1]
\State {\bf Input:} a set of random terminal functions $G:=\{g^k(x)\}_{k=1}^{N_{\text{train}}}$, number of policy iterations $M$, random sample points $x_1,...,x_{N_0}\in\mathbb{R}^d$, initial control $u_0$, weight of loss $\alpha_1$ and $\alpha_2$, arbitrary terminal function $g^{\text{new}}(x)$;
%\State $t \leftarrow 1$, $x_1 \leftarrow 0$, $\mathcal{D} \leftarrow \emptyset$;
\State {\bf Output:} operator $\mathcal{H}_\theta(g^{\text{new}}(x))$ solving~\eqref{eq:app_hj} with $g(x)=g^{\text{new}}(x)$
\For{$n=0,1,2,\ldots,M-1$} 
\While{trained}
\State With all $g^k \in G$, $V_n^h:=\mathcal{G}_\theta [g^k,u_n](t,x)$.
\State $\mathcal{L}_1 = (\sum_{i,j,k}\partial_t V_n^h+\nabla^h V_n^h \cdot f(t_i,x_j,u_n)+Nh\Delta^h V_n^h)^2$.
\State $\mathcal{L}_2 = \sum_{j,k}(V_n^h(T,x_j)-g^k(x_j))^2$.
\State $\theta \leftarrow \text{Adam}(-\alpha_1\mathcal{L}_1-\alpha_2 \mathcal{L}_2)$
\EndWhile
\State $V_n^h \leftarrow \mathcal{G}_\theta[g^k,u_n]$
\State (Update control)
\[
    u_{n+1}(t,x)=\argmin_{{\boldsymbol{u}}\in U} \{\nabla^h V_n^h \cdot f(t,x,{\boldsymbol{u}})+L(t,x,{\boldsymbol{u}})\},
\]
\EndFor
\State Return $\mathcal{H}_\theta[g^{\text{new}}]:=V_M^h$
\end{algorithmic}
\end{algorithm}

To find the optimal control through the lens of policy iteration, we iteratively employ the PI-DeepONet as illustrated in Figure \ref{fig:framework}. The left side of the figure illustrates a single iteration step for approximating the solution of \eqref{eq:semidiscret} using PI-DeepONet $\mathcal{G}_\theta$, while the right side represents the process of repeating this process iteratively. 
\begin{equation}\label{eq:twostep}
    \begin{aligned}
    V_n^h(t,x) &\approx \mathcal{G}_\theta[g,u_n],&&\text{(via PI-DeepONet)}, \\
    u_{n+1}(t,x)&=\argmin_{{\boldsymbol{u}}\in U} \{\nabla^h V_n^h \cdot f(t,x,{\boldsymbol{u}})+L(t,x,{\boldsymbol{u}})\},&&\text{(via optimization)},
\end{aligned}
\end{equation}
for $n=0,1,...M-1$ with the initial value $u_0(t,x)$. Here, $u_0(t,x)$ is set to be 0 for convenience of analysis. After $M$ iteration steps, we end up with $V_M^h$ and define 
\[
\mathcal{H}_\theta(g):=V_M^h.
\]

Owing to the theoretical result of the policy iteration as discussed in the previous section, we deduce that $V_M^h(t,x)$ approximates discrete-in-space viscosity solution $V^h$ for sufficiently large $M$. As a byproduct, the optimal control $u_*$ associated with $V$ is achieved with $V_M^h$.  
Note that we use $N_\text{train}$ randomly chosen terminal functions that are given by $\{g^i(x)\}_{i=1}^{N_\text{train}}$ to train the PI-DeepONet, and exploit the trained PI-DeepONet structure to infer the solution to the optimal control problem corresponding to a new terminal function $g^{\text{new}}(x)$. Empirical evidence for the effectiveness of this approach is provided in Section~\ref{sec:exp}. The detailed procedure is found in Algorithm~\ref{alg:main}.

Let us discuss the potential advantages of our new approach. (i) The Hamiltonian $H$ is computed only at some sample grid points. (ii) The implementation of an operator learning framework allows us to fast inferences after the model is trained. (iii) The temporal derivative can be computed effectively via an automatic differentiation tool, Autograd from Pytorch~\cite{paszke2017automatic}. (iv) The spatial derivative is not required as it is replaced by the finite difference.

%Therefore, our framework approximates the operator $\mathcal{G}:g\mapsto\mathcal{G}[g]=V^h$ by iteratively leveraging the PI-DeepONet $\mathcal{G}_\theta:[g(x),u_n(t,x)]\mapsto V_n^h(t,x)$. We denote the mapping from $g$ to $V^h_M$ by $\mathcal{H}_\theta[g]:=V^h_M$.

\subsection{Error analysis}
In this subsection, we present an error analysis for the solution obtained via our algorithm. The PI-DeepONet framework trains the operator in a way that minimizes both PDE loss and terminal condition loss. Assuming the trained error is small, one can examine the convergence of the solution trained via our main algorithm to the actual viscosity solution $V$ of~\eqref{eq:hj}.

Throughout the section, we impose the following assumption which indicates that both PDE and boundary loss diminish as the training proceeds. 

\begin{enumerate}[label=(A3)]
\item For any $g$ satisfying (A2) and $n \geq 0$, let $v_n(t,x;\theta):=\mathcal{G}_\theta[g,u_n]\in C^1([0,T]\times \mathbb{R}^d)$ be obtained from~\eqref{eq:twostep} and satisfy in the classical sense that
\[
\begin{cases}
\|\partial_t v^h_n+ L(t,x,u_n)+\nabla^h v_n^h \cdot f(t,x,u_n)+N h \Delta^h v_n^h \|_\infty \leq \epsilon_{1,n},\\
\|g(x)-v_n^h(T,x)\|_\infty \leq  \epsilon_{2,n}
\end{cases}
\]
for some $\epsilon_{1,n} \searrow 0$ and $\epsilon_{2,n}\searrow 0$ and
\[
u_{n+1}(t,x)=\argmin_{{\boldsymbol{u}}\in U}\{\nabla^h v^h_n \cdot f(t,x,{\boldsymbol{u}})+L(t,x,{\boldsymbol{u}})\}.
\]
\end{enumerate}

By the equivalence of the viscosity solution and classical solution for continuously differentiable solutions, we further have that $v_n^h$ satisfies 
\begin{equation}\label{eq:super}
\begin{cases}
\partial_t v^h_n+ L(t,x,u_n)+\nabla^h v_n^h \cdot f(t,x,u_n)+N h \Delta^h v_n^h  \leq \epsilon_{1,n},\\
v_n^h(T,x)-g(x) \geq  -\epsilon_{2,n},
\end{cases}
\end{equation}
and 
\begin{equation}\label{eq:sub}
\begin{cases}
\partial_t v^h_n+ L(t,x,u_n)+\nabla^h v_n^h \cdot f(t,x,u_n)+N h \Delta^h v_n^h  \geq -\epsilon_{1,n},\\
v_n^h(T,x) - g(x)\leq \epsilon_{2,n},
\end{cases}
\end{equation}
in the viscosity sense.

Based on this observation, we slightly modify the argument in~\cite{tang2023policy} to verify the stability of $v_n^h$ satisfying both equations~\eqref{eq:super} and \eqref{eq:sub} in viscosity sense as $n$ grows to infinity. 
\begin{prop}[Uniform boundedness]\label{prop:bound}
Let $T>0$ be given and assume (A1)-(A3). Let $N \geq \max\{1,\|f\|_\infty/2\}$, we deduce that $v_n^h \in C^1([0,T]\times \mathbb{R}^d)$ are uniformly bounded for all $n\geq 0$ and $h>0$.
\end{prop}
\begin{proof}
It is enough to show that $v_0^h$ is uniformly bounded as $v_n^h$ is defined iteratively. Since $\|L\|_\infty,\|g\|_\infty <\infty$, 
$v_0^h \pm \epsilon_{1,n}(T-t)\pm\epsilon_{2,n}$ are supersolution and subsolution to 
\[
\begin{cases}
\partial_t v + L(t,x,u_0)+\nabla^h v \cdot f(t,x,u_0)+Nh\Delta^h v =0 \quad&\text{in}\quad (0,T)\times \mathbb{R}^d,\\
v(T,x)=g(x)\quad&\text{on}\quad \mathbb{R}^d,
\end{cases}
\]
whose unique viscosity solution is denoted by $v$. By the comparison principle,
\[
v_0^h - \epsilon_{1,n}(T-t) - \epsilon_{2,n}\leq v \leq v_0^h + \epsilon_{1,n}(T-t) + \epsilon_{2,n},
\]
and combining with the uniform boundedness of $v$~\cite{tang2023policy}[Proposition 2.2], we derive that $v_0^h$ is uniformly bounded. Applying the argument iteratively, we deduce that $v_n^h$ is uniformly bounded for each $n$ and $h$.
\end{proof}

We define the following quantity to quantify the cumulative error arising from~\eqref{eq:super} and~\eqref{eq:sub}.

\begin{defn}\label{def:ep_t}
With $\{\epsilon_{i,n}\}$ for $i=1,2$ satisfying (A3) and $T>0$ fixed, let us define
\[
\epsilon(t) := \lim_{n\rightarrow \infty}\{(T-t)(\epsilon_{1,0} + 2(\sum_{m=2}^{n-1} \epsilon_{1,m}) +\epsilon_n)+(\epsilon_{2,0} + 2(\sum_{m=2}^{n-1} \epsilon_{2,m}) +\epsilon_{2,n})\},
\]
where $t\in [0,T]$.
\end{defn}
The term $\epsilon(t)$ is determined by the cumulative loss that occurred in the training of the algorithm and plays a central role in the error analysis.
\begin{thm}[Stability of $v_n^h$]\label{prop:stability}
Let $\epsilon(t)<\infty$ be given as in Definition~\ref{def:ep_t}. Under the same assumptions as in Proposition~\ref{prop:bound}, 
$\{v_n^h(\cdot;\theta)\}$ converges locally uniformly as $n$ increases. 
\end{thm}

\begin{proof}
For $u\in \mathcal{U}_t$, define an operator
\[
\mathcal{L}^h(v,u):=\partial_t v +c(t,x,u)+\nabla^h v \cdot f(t,x,u)+N h \nabla^h v.
\]
By the definition of $u_n$, we have that
\[
c(t,x,u_{n+1})+\nabla^h v_n^h \cdot f(t,x,u_{n+1})\leq c(t,x,u_{n})+\nabla^h v_n^h \cdot f(t,x,u_{n}).
\]
Hence,
\[
\mathcal{L}^h[v_{n}^h,u_{n+1}] \leq \epsilon_{1,n}.
\]
Clearly, $v_n^h+\epsilon_{2,n}+\epsilon_{2,n+1}+(\epsilon_{1,n}+\epsilon_{1,n+1})(T-t)$ and $v_{n+1}^h$ are supersolution and subsolution to
\begin{equation*}
\begin{cases}
\mathcal{L}[v,u_{n+1}] = -\epsilon_{1,n+1},\\
v(T,x)-g(x)=\epsilon_{2,n+1}.
\end{cases}
\end{equation*}
Again by the comparison principle,
\[
v_{n+1}^h \leq v_n^h+(\epsilon_{2,n}+\epsilon_{2,n+1})+(\epsilon_{1,n}+\epsilon_{1,n+1}) (T-t),
\]
which can be written as
\[
(v_{n+1}^h-\beta_{n+1})- (v_n^h -\beta_n) \leq 0,
\]
for
\[
\beta_n(t) := (T-t)(\epsilon_{1,0} + 2(\sum_{i=2}^{n-1} \epsilon_{1,i}) +\epsilon_{1,n})+(\epsilon_{2,0} + 2(\sum_{i=2}^{n-1} \epsilon_{2,i}) +\epsilon_{2,n}).
\]
Therefore, by the uniform boundedness and monotonicity of $\{(v_n^h(t,x) -\beta_n(t)\}_{n=0}^\infty$,  
\[
\lim_{n\rightarrow \infty} v_n^h(t,x) - \beta_n(t) := \tilde v^h(t,x) \quad \text{locally uniformly}.
\]
Hence, we obtain the following locally uniform convergence of $v_n^h$, which is
\[
\lim_{n\rightarrow \infty} v_n^h(t,x) := \tilde v^h(t,x) + \epsilon(t),
\]
where
\[
\epsilon(t) = \lim_{n\rightarrow \infty}\{(T-t)(\epsilon_{1,0} + 2(\sum_{i=2}^{n-1} \epsilon_{1,i}) +\epsilon_n)+(\epsilon_{2,0} + 2(\sum_{i=2}^{n-1} \epsilon_{2,i}) +\epsilon_{2,n})\}.
\]
\end{proof}

Finally, one can further verify the following stability result on the approximate solution obtained from our algorithm.
%\[
%\lim_{n\rightarrow \infty} v_n^h(t,x;\theta) \quad\text{exists for any }\theta.
%\]

\begin{cor}
Under assumptions (A1)-(A3) and $N\geq \|f\|_\infty/2$, $h \in (0,1)$, let $V$ be a unique viscosity solution to~\eqref{eq:hj} and 
\[
v^h(\cdot;\theta)=\lim_{n\rightarrow \infty} v_n^h(\cdot;\theta).
\]
Then we have that
\[
\|V-v^h\|_\infty \leq C\sqrt{h}.
\]
\end{cor}
\begin{proof}
%Let $\epsilon(t)$ and $\tilde v^h(t)$ be from the proof of Theorem~\ref{prop:stability}, that is,
%\[
%\lim_{n\rightarrow \infty} v_n^h(t,x) := \tilde v^h(t,x) + \epsilon(t)
%\]
%for $\epsilon(t)<\infty$. 
By the stability property of viscosity solution~\cite{achdou2013introduction} and the uniqueness of a viscosity solution to~\eqref{eq:app_hj}, we have that $V^h \equiv v^h(t,x;\theta)$, which implies that $v^h$ is independent of the neural network parameter $\theta$.

Therefore, we have that
\begin{equation*}
\begin{split}
|V-v^h | &\leq |V-V^h|+\lim_{n\rightarrow \infty}|V^h - v_n^h|\\
&\leq C\sqrt{h} 
\end{split}
%+ \lim_{n\rightarrow\infty}|V^h - %v_n^h- \epsilon(t)|+\epsilon(t)\\
%& \leq C\sqrt{h}+ |V^h - \tilde v^h | +\epsilon(t)\\
%& %\leq C\sqrt{h} + \epsilon(0)
\end{equation*}
for some $C>0$ where the second inequality comes from~\cite{tang2023policy}[Theorem 3.3]. %and the third inequality follows from the fact that $\lim_{n\rightarrow\infty} v_n^h+\epsilon(t) \equiv V_h$ thanks to the uniqueness of viscosity solution to~\eqref{eq:app_hj}.
\end{proof}

\begin{comment}
The following states that well-trained operator yields a solution close to the actual solution $v$.
At each intermediate step, we approximate $v^h$ by the neural network solution $v_N$ minimizing the loss.
\end{comment}

\begin{comment}
\begin{thm}
If trained well, there exists a subsequence $n_k$ such that
\[
\lim_{k\rightarrow\infty} v^{n_k}(\cdot;\theta) =v,
\]
and hence,
\[
\lim_{k\rightarrow \infty} |v^{n_k}(\cdot;\theta)- V|
\]
\end{thm}
\end{comment}

\section{Experimental result}\label{sec:exp}
We provide results of an extensive set of experiments that verify the effectiveness of the proposed method. Our framework offers multiple advantages. We first emphasize that once the neural network is trained, not only the value function $V(t,x)$ that is a solution to the corresponding HJB equation but also an optimal trajectory-control pair $(x(s),u(s))$ is obtained systematically for any terminal function $g$ given. Additionally, we avoid the curse of dimensionality since only a few measurements are used instead of discretizing the whole high-dimensional space, which is one of the well-known advantages of exploiting the framework of PI-DeepONet as highlighted in~\cite{goswami2022physics}.

% \subsection{Commodity trading}
% \begin{itemize}
%     \item $x_1(t)$ = money on hand at time t
%     \item $x_2(t)$ = amount of wheat owned at time t 
%     \item $\alpha(t)$= rate of buying or selling of wheat 
%     \item $q(t)$ = price of wheat at time t (known)
%     \item $\lambda$ = cost of storing a unit amount of wheat for a unit of time.
% \end{itemize}
% $|\alpha(t)| \leq M$. Cost is defined as 
% \[
% J[x,\alpha]:= x_1(T) +q(T)x_2(T)
% \]
% with $[x_1(0);x_2(0)]= [x_1;x_2]$.
% The dynamics is given by
% \[
% \dot x_1 =-\lambda x_2-q(t)\alpha(t), \quad \dot x_2 =\alpha
% \]

% Optimal control is known to have a form
% \begin{equation}
%     \alpha(t)=
%     \begin{cases}
%     M\quad\text{if}\quad q(t)<p^2(t)\\
%     -M \quad\text{if}\quad q(t)>p^2 (t)
%     \end{cases}
% \end{equation}
% for $p^2(t):= \lambda (t-T)+q(T)$

% \subsection{(i) Synthetic case 1D} \cite{wang2003numerical}
% \[
% v(x,s)=\min_{\mathcal{U}} \{-x^2(1)\}
% \]
% subject to
% \begin{align*}
% \frac{dx}{dt}= u(t)
% \end{align*}
% where $\mathcal{U}=\{u:[0,1] \rightarrow [0,1]\}$.
% Then the true solution is given by
% \begin{equation*}
% v(t,x)=
% -[|x|-(1-t)]^2
% \end{equation*}

% \subsection{(ii) Synthetic case 1D}
% \[
% J[x,\alpha]=-x(2)
% \]
% Dynamics is given by
% \[
% \dot x = \frac{5}{2}(-x^2 + ux -u^2)
% \]
% with $x(0)=1$.
% Exact trajectory is given by
% \[
% x(t)=4/(1+3e^{5t/2})
% \]
% and control is
% \[
% u(t)=x(t)/2
% \]

% \subsection{Synthetic case 5D}

\subsection{Synthetic case 2D}
Let us introduce a 2D nonlinear vehicle example that has been considered in \cite{lee2021computationally}. In this example, the dynamics of a vehicle is given by
\begin{equation}\label{vehic:dy}
\begin{cases}
\frac{d x_t }{d t} = \begin{bmatrix}
    \cos (u_t) \\ \sin (u_t)
\end{bmatrix}, \quad t\in (0,T),\\
x_t=x,\\
\end{cases}
\end{equation}
where $x_t \in \mathbb{R}^2$ and $u_t \in U=[-\pi,\pi]$ denote trajectory and control respectively.  Our goal is to understand the following optimization problem:
\begin{equation}\label{eq:vehicle}
V(t,x)= \inf_{u \in \mathcal{U}_t} \|x_T \|_2,
\end{equation}
for $\mathcal{U}_t= \{u:[t,T]\mapsto [-\pi,\pi] : u \text{ is measurable}\}$ subject to the dynamics above. This optimal control problem aims at determining an optimal steering angle at each time $t$ that drives the vehicle moving with speed 1 to the origin as close as possible. Recalling~\eqref{eq:update}with 
\[
f(t,x,\boldsymbol{u})=\begin{bmatrix}
    \cos(\boldsymbol{u}) \\ \sin(\boldsymbol{u})
\end{bmatrix}\quad\text{and}\quad L(t,x,\boldsymbol{u})=0,
\]
we have the explicit formula for $u_{n+1}$:
\begin{align*}
u_{n+1}(t,x)&=u(t,x,\nabla_{x}^hv_n^h(t,x))\\
&=\begin{cases}
    p-\pi&\text{if}\;\nabla_{1}^hv_n^h>0,\nabla_{2}^hv_n^h>0,\\
    p+\pi&\text{if}\;\nabla_{1}^hv_n^h>0,\nabla_{2}^hv_n^h<0,\\
    p &\text{else},
\end{cases}
\end{align*}
where $p=\arctan\left(\frac{\nabla_{2}^hv_n^h(t,x)}{\nabla_{1}^h v_n^h(t,x)}\right)$ ($-\frac{\pi}{2}< p < \frac{\pi}{2}$).

For the experiment, we use $T=1$. Since the velocity is set to be 1, it is clear that the optimal control is given by $u_t=\arctan\big(\frac{x[2]}{x[1]}\big)$ if $\|x\|_2>1$ and $x[1]<0$, which means that the vehicle moves toward the origin from the initial point $x_0$. Our goal is to find an approximate solution $V_h$ to HJB equation
\begin{equation*}
\begin{cases}
\partial_t V + \inf_{  \boldsymbol{u}\in U } \bigg\{ \begin{bmatrix} \cos(u)\\ \sin(u)\end{bmatrix} \cdot D_x V \bigg\}=0 \quad &\text{in}\quad (0,1)\times \mathbb{R}^2,\\
V(1,x)=\|x\|^2_2 \quad&\text{on}\quad \mathbb{R}^2,
\end{cases}
\end{equation*}
and synthesize optimal control accordingly. To compute $V(0,x)$ numerically, we discretize the time to solve~\eqref{eq:vehicle} subject to~\eqref{vehic:dy}, that is,
\[
V(0,x)\approx \inf \|x_1\|_2
\]
subject to
\[
x_{t_{k+1}}=x_{t_k} +(t_{k+1}-t_k) \begin{bmatrix}
    \cos(u_{t_k})\\\sin(u_{t_k})
\end{bmatrix} \quad \text{and} \quad x_{0}=x
\]
for $0=t_0<...<t_\ell=1$ with $\ell \in \mathbb{N}$. The interior point method~\cite{nocedal1999numerical} is used to solve for optimization variables $(x_{t_k},u_{t_k})$ for $k \in [0,\ell-1]$. For the rest of the parameters, we set $h=0.005$, $M=5$, $N=1$. 

\begin{comment}
\begin{equation}
g(x_1,x_2)=2+a_1\cos(2\pi(b_1x_1+b_2x_2))+a_2\sin(2\pi(b_3x_1+b_4x_2)),    
\end{equation}
\end{comment}

\begin{figure}[t!]
    \centering
    \includegraphics[width=\textwidth]{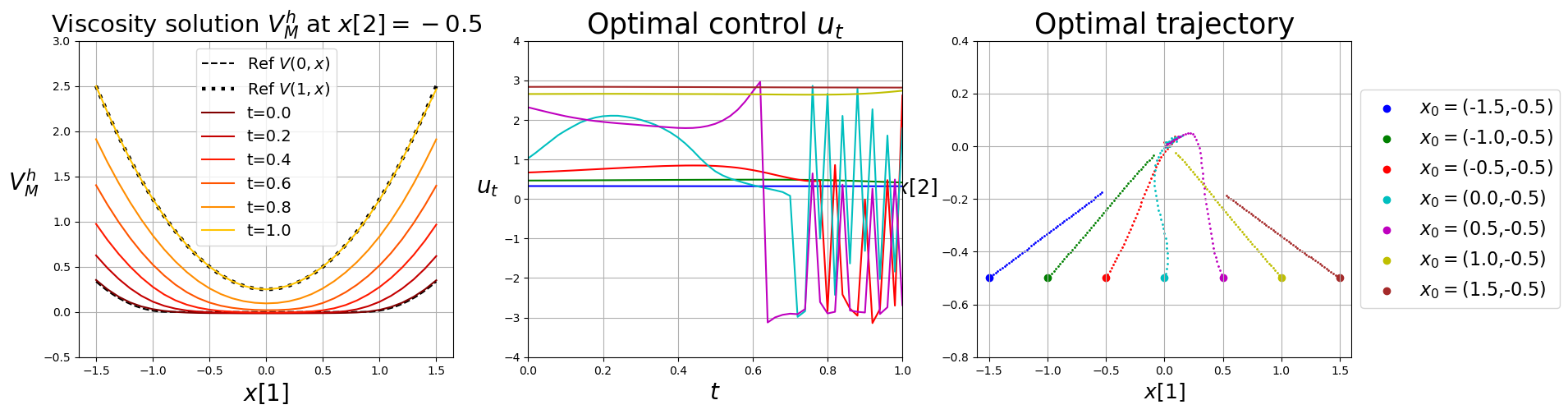}
    \caption{Plots of $V_M^h(t,x)$ with $t \in [0,1]$ fixing $x[2]=-0.5$ (left). Optimal controls corresponding to different initial points of $x$ with $x[1] \in [-1.5,1.5]$ with $x[2]=-0.5$ (middle). Optimal state trajectories associated with the controls computed (right).}
    \label{fig:vehicle}
\end{figure}

The left plot in Figure~\ref{fig:vehicle} demonstrates that the approximated viscosity solution $V_M^h$ closely matches the reference numerical solution. The middle and right plots in Figure~\ref{fig:vehicle} depict the optimal control $u_t$ and the optimal trajectory for different initial points $x_0$. In the case of $x_0 = (-1.5, -0.5)$ and $x_0 = (-1, -0.5)$, we expect the optimal actions to be $\arctan(1/3)$ and $\arctan(1/2)$, respectively heading to the origin. The middle plot also verifies such a choice of steering angles. On the other hand, when $x_0 = (0.5, -0.5)$ and $x_0 = (0, -0.5)$, the vehicle has sufficient time to reach the origin. For these initial points, we observe oscillations in control as shown in the right plot in Figure~\ref{fig:vehicle}.

\subsection{LQR with a compact control set}
We consider a LQR problem given as
\begin{equation}\label{eq:lqr}
\begin{cases}
\frac{dx_t}{dt}&= Ax_t +B u_t, \quad  t \in (0,T),\\
x_t&=x,
\end{cases}
\end{equation}
where $x_t\in \mathbb{R}^d$ and $u_t \in U$ for some compact set $U\subset \mathbb{R}^m$ denote the state and control respectively. The value function is defined as 
\begin{equation}\label{eq:lqr_exp}
V(t,x):=\inf_{u \in \mathcal{U}_t} \left \{  \int_t^T (x_s^\top Q x_s+u_s^\top R u_s ) \de s+g(x_T) : x_t = x \right \},
\end{equation}
where $g$ is a given function and $Q,R$ are positive definite. The Hamiltonian is given by
\[
H(t, x, p):=\inf_{{\boldsymbol{u}}\in U} \left\{ p\cdot f(t,x,{\boldsymbol{u}})+L(t,x,\boldsymbol{u}) \right\},
\]
with $f(t,x,{\boldsymbol{u}})=Ax+B{\boldsymbol{u}}$ and $L(t,x,\boldsymbol{u})=x^\top Q x + \boldsymbol{u}^\top R \boldsymbol{u}$. 

\noindent\textbf{Experimental setup.} We consider two cases; $(d,m)=(5,3)$ and $(d,m)=(10,5)$ with $U=[-1/3,1/2]^m$. For system matrices $A$ and $B$, each component of $A$ and $B$ are randomly chosen from $[0,0.2]$ and $[0,0.1]$ respectively. We use the identity matrices for $Q$ and $R$. For $g(x)$, we take $g_k(x)=0.3+0.1 \times k\|x\|_2^2$ with $k=1,2,3$ to let the model be trained with various terminal functions. We emphasize that the solution to the HJB equation can be inferred with different terminal functions leading to different optimal control trajectories. For other parameters, we use $T=0.5$, $h=0.005$, $N=1$, and $M=3$ regardless of the dimension.

\noindent\textbf{Comparison with the true solution.} Similar to the toy example, we provide the comparison between $V_M^h(0,x)$ and $V(0,x)$. To compute $V(0,x)$, we discretize the time to solve~\eqref{eq:lqr_exp} subject to~\eqref{eq:lqr}, that is,
\[
V(0,x)\approx \inf_{x_{t_k},u_{t_k}}\{\sum_{k=0}^{\ell-1} x_{t_k}^\top Q x_{t_k} + u_{t_k}^\top R u_{t_k} +g(x_{0.5}) \}
\]
subject to
\[
x_{t_{k+1}}=x_{t_k} +(t_{k+1}-t_k)(A x_{t_k} +B u_{t_k}) \quad \text{and} \quad x_{0}=x
\]
for $0=t_0<...<t_\ell=0.5$ with $\ell \in \mathbb{N}$. Again we use the interior point method~\cite{nocedal1999numerical} to solve for optimization variables $(x_{t_k},u_{t_k})$ for $k \in [0,\ell-1]$.

\subsubsection{$(d,m)=(5,3)$ with $A \in \mathbb{R}^{5\times 5}$ and $B \in \mathbb{R}^{5\times 3}$.}
For the dynamics, we use
\[\footnotesize
A=
  \begin{bmatrix}
    0.08 & 0.01 & 0.01 & 0.15 & 0.05 \\
    0.16 & 0.16 & 0.15 & 0.18 & 0.16\\
    0.14 & 0.17 & 0.07 & 0.08 & 0.12\\
    0.12 & 0.15 & 0.11 & 0.18 & 0.02\\
    0.12 & 0.08 & 0.13 & 0.1  & 0.09\\
  \end{bmatrix}
\]
and
\[\footnotesize
B=
\begin{bmatrix}
0  & 0.05 & 0.06\\
0.07 & 0.01 & 0.04\\
0.02 & 0.   & 0.1 \\
0.09 & 0.08 & 0.08 \\
0.01 & 0.07 & 0.05 \\
\end{bmatrix}.
\]

\begin{figure}[t!]
    \centering
    \includegraphics[width=\textwidth]{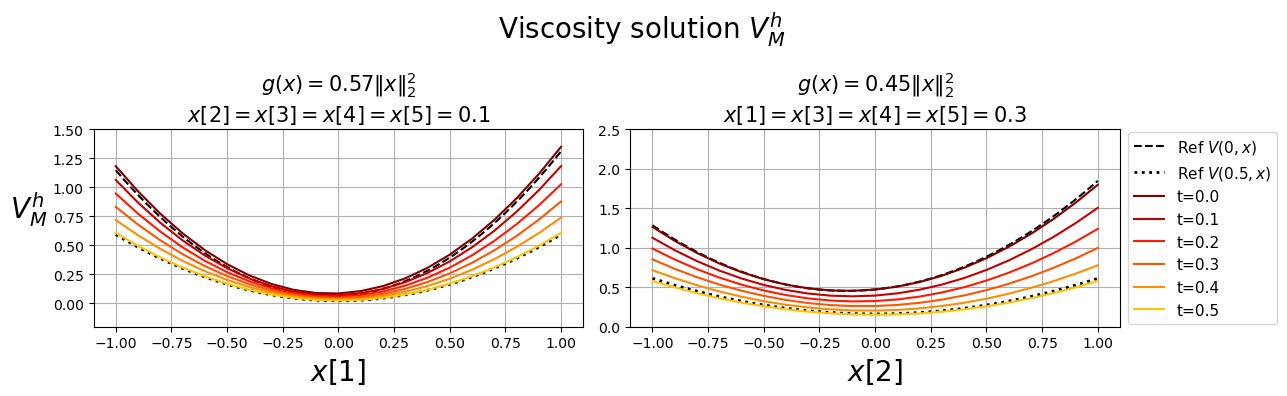}
    \includegraphics[width=\textwidth]{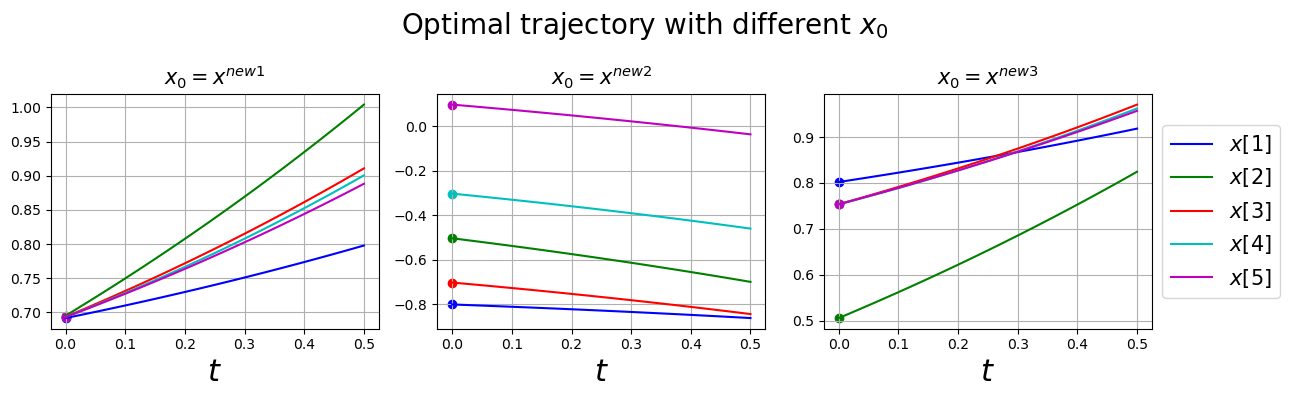}
    \includegraphics[width=\textwidth]{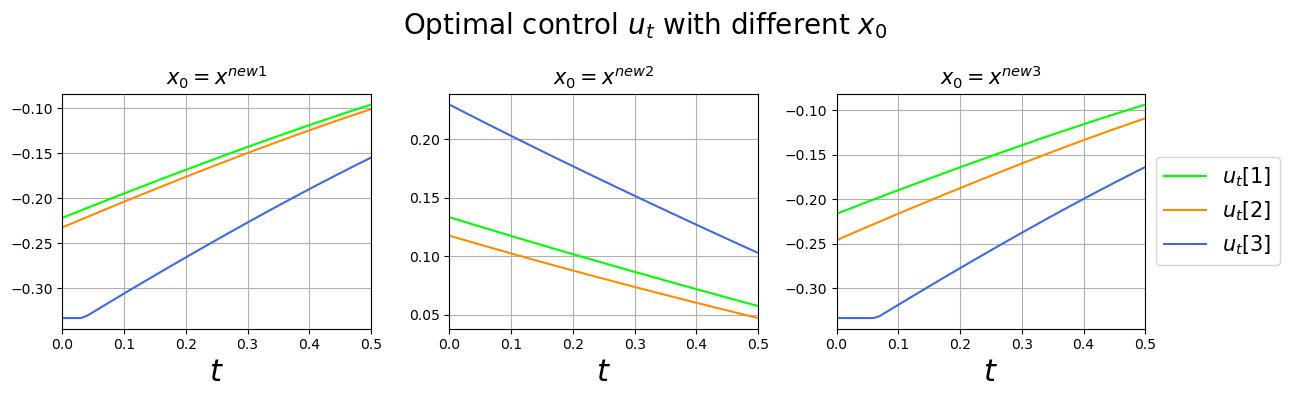}
    \caption{Viscoisity solutions, state, and control trajectories for LQR problem with $(d,m)=(5,3)$. Value functions $V_M^h$ with different terminal functions $g(x)$ for the LQR problem (top). Optimal trajectories (middle) and controls (bottom) associated with different initial values of $x_0= x^{\text{new}}$.}
    \label{fig:lqr_5_3_g06}
\end{figure}

The results are presented in Figure~\ref{fig:lqr_5_3_g06}. The evolution of the value function $V_M^h(t,x)$ as $t$ and a component of $x$ vary, while fixing the rest $x$, is shown in the first row of Figure~\ref{fig:lqr_5_3_g06}. It shows that the approximated value function $V_M^h(t,x)$ aligns with $V(t,x)$ at $t=0$ and $t=0.5$. To verify that our learned operator efficiently approximates the value function for unseen inputs, we use the terminal functions $g(x)=0.57 |x|_2^2$ and $g(x)=0.45 |x|_2^2$. The second and third rows of Figure~\ref{fig:lqr_5_3_g06} show the optimal trajectory and control with three different initial points $x_0$.

%\marginpar{For different $g(x)$s, (1) plot viscosity solution (2) plot optimal control for each component $u\in\mathbb{R}^5$ for several random points.} 

\subsubsection{$(d,m)=(10,5)$ with $A \in \mathbb{R}^{10\times 5}$ and $B \in \mathbb{R}^{5\times 5}$.}
Let us take
\[\tiny
A=
\begin{bmatrix}
0.00 & 0.15 & 0.01 & 0.11 & 0.03 & 0.19 & 0.06 & 0.07 & 0.14 & 0.07 \\
0.06 & 0.10 & 0.14 & 0.09 & 0.16 & 0.01 & 0.15 & 0.17 & 0.09 & 0.11 \\
0.11 & 0.02 & 0.10 & 0.19 & 0.04 & 0.14 & 0.18 & 0.01 & 0.10 & 0.16 \\
0.15 & 0.13 & 0.01 & 0.17 & 0.04 & 0.06 & 0.16 & 0.03 & 0.08 & 0.15 \\
0.02 & 0.19 & 0.19 & 0.17 & 0.13 & 0.15 & 0.00 & 0.17 & 0.08 & 0.17 \\
0.07 & 0.00 & 0.02 & 0.14 & 0.10 & 0.08 & 0.13 & 0.07 & 0.03 & 0.05 \\
0.00 & 0.15 & 0.16 & 0.12 & 0.17 & 0.06 & 0.05 & 0.14 & 0.18 & 0.10 \\
0.02 & 0.10 & 0.14 & 0.12 & 0.17 & 0.01 & 0.15 & 0.08 & 0.10 & 0.17 \\
0.04 & 0.03 & 0.07 & 0.02 & 0.13 & 0.10 & 0.01 & 0.13 & 0.15 & 0.09 \\
0.00 & 0.12 & 0.07 & 0.01 & 0.09 & 0.15 & 0.06 & 0.05 & 0.08 & 0.05
\end{bmatrix}
\]
and
\[\tiny
B=
\begin{bmatrix}
0.02 & 0.05 & 0.09 & 0.08 & 0.06 \\
0.05 & 0.06 & 0.10 & 0.10 & 0.01 \\
0.06 & 0.01 & 0.05 & 0.04 & 0.05 \\
0.02 & 0.09 & 0.00 & 0.03 & 0.07 \\
0.07 & 0.09 & 0.02 & 0.05 & 0.05 \\
0.00 & 0.02 & 0.03 & 0.05 & 0.03 \\
0.04 & 0.07 & 0.03 & 0.01 & 0.03 \\
0.10 & 0.05 & 0.02 & 0.05 & 0.03 \\
0.08 & 0.06 & 0.08 & 0.06 & 0.09 \\
0.00 & 0.00 & 0.06 & 0.04 & 0.05
\end{bmatrix}
\]

\begin{figure}[t!]
    \centering
    \includegraphics[width=\textwidth]{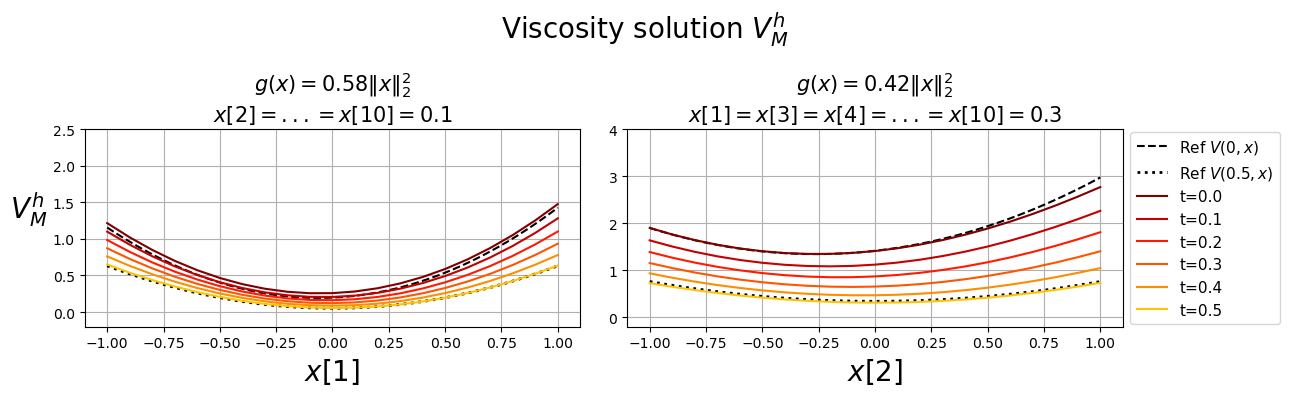}
    \includegraphics[width=\textwidth]{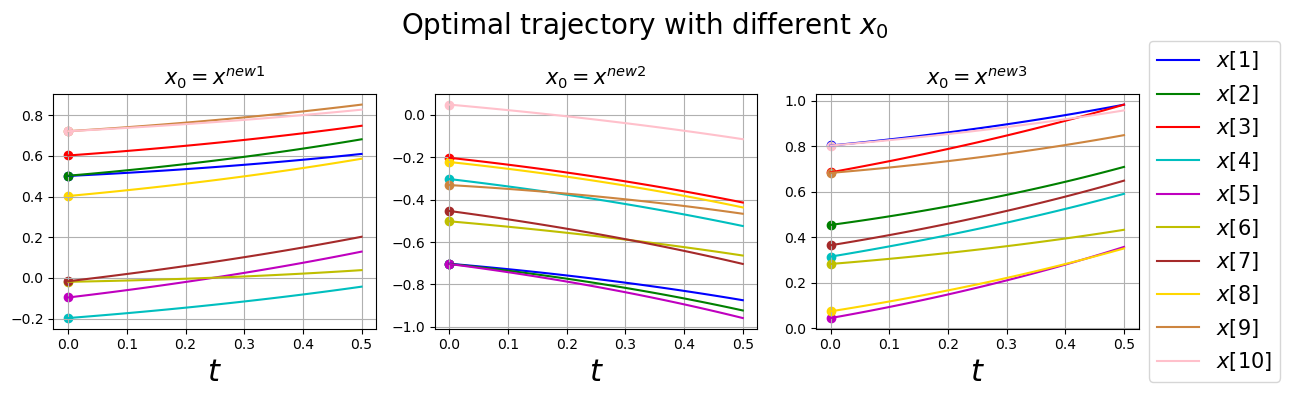}
    \includegraphics[width=\textwidth]{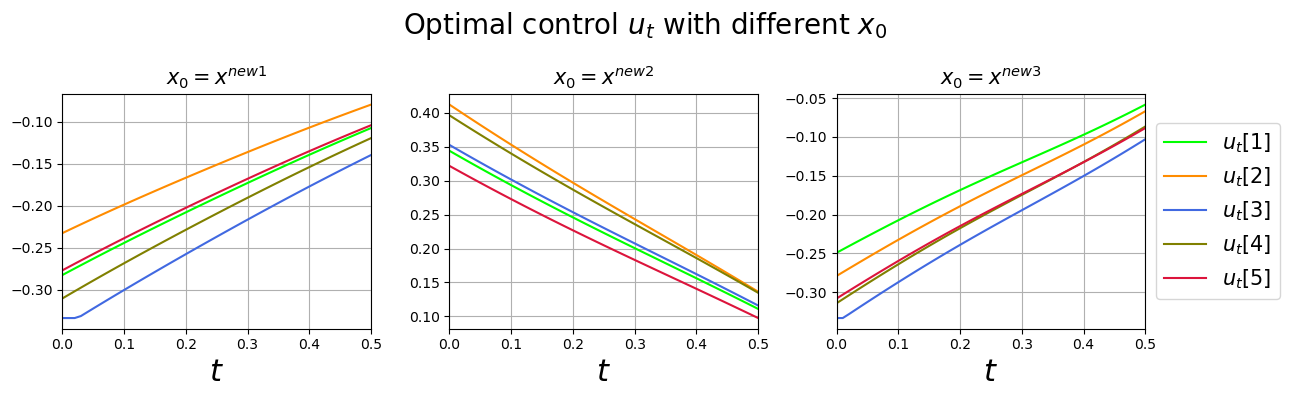}
    \caption{Viscoisity solutions, state, and control trajectories for LQR problem with $(d,m)=(5,3)$. Value functions $V_M^h$ with different terminal functions $g(x)$ for the LQR problem (top). Optimal trajectories (middle) and controls (bottom) associated with different initial values of $x_0= x^{\text{new}}$.}
    \label{fig:lqr_10_5_g06}
\end{figure}

We use slightly different terminal functions $g(x)=0.58\|x\|_2^2$ and $g(x)=0.42\|x\|_2^2$ to verify the robustness of inferring the solution $V_M^h$ associated with an unseen terminal function $g(x)$. As found in first row of Figure~\ref{fig:lqr_10_5_g06}, the learned value function $V_M^h(t,x)$ matches with $V(t,x)$ for $t=0$ and $t=0.5$. Optimal trajectories and controls are computed for three different initial points as shown in the second and third rows of Figure~\ref{fig:lqr_10_5_g06}.

% For the optimization, we contatenate controls $\{\alpha_i\}_{i=1}^B$ where $B$ denotes the batch size to construct
% {
% \color{red}
% \[
% {u} := \begin{bmatrix}  \alpha_1^\top & \cdots & \alpha_B^\top\end{bmatrix}^\top,
% \]
% }
% and perform batch-wise optimization for acceleration. 
\section{Conclusion}\label{sec:conc}
We propose a framework for solving HJB equations and optimal control problems with the aid of PI-DeepONet. By leveraging the idea of policy iteration with semi-discretization whose convergence result is known, PI-DeepONet can be applied successfully iterative manner. The uniform error bound suggests that the algorithm is stable. We emphasize that this approach is capable of solving multiple HJB equations and optimal control problems. Various experimental results support the effectiveness of our algorithm. For future work, we propose to explore the infinite horizon problem and static HJ equations where such an iteration has not been explored.

\section*{Acknowledgement}
The authors are thankful to Hung Vinh Tran at the University of Wisconsin-Madison for sharing ideas and discussions. This work was supported by the National Research Foundation of Korea (NRF) grant funded by the Korea government(MSIT) (RS-2023-00219980,RS-2023-00211503) Jae Yong Lee was supported by Institute for Information \& Communications Technology Planning \& Evaluation (IITP) through the Korea government
(MSIT) under Grant No. 2021-0-01341 (Artificial Intelligence Graduate School Program (Chung-Ang University)).

 \bibliographystyle{elsarticle-num} 
 \bibliography{ref}





% \end{thebibliography}
\end{document}